\documentclass[12pt,twoside]{article}
\usepackage{amsmath, amssymb, amsthm}
\font\fourteenb=cmb10 at 14pt
\begin{document}

\vspace*{-1.0cm}\noindent \copyright
 Journal of Technical University at Plovdiv\\[-0.0mm]\
\ Fundamental Sciences and Applications, Vol. 11, 2005-2006\\[-0.0mm]
\textit{Series A-Pure and Applied Mathematics}\\[-0.0mm]
\ Bulgaria, ISSN 1310-8271\\[+1.2cm]
\font\fourteenb=cmb10 at 12pt
\begin{center}

{\bf \LARGE Construction of analytic functions, which determine
bounded
Toeplitz operators on  ${\bf H}^{{\bf 1}} $  and  ${\bf
H}^{{\bf \infty }} $
\\ \ \\ \large Peyo Stoilov}
\end{center}

\footnotetext{{\bf 1991 Mathematics Subject
Classification:} Primary 30E20, 30D50} \footnotetext{{\it Key words
and phrases:} Analytic function, Toeplitz operators, Cauchy integrals, multipliers.
}\footnotetext{{\it Received June 15, 2005.}}
\begin{abstract}
For  $f\in H^{\infty } $  we denote by  $T_{f} $  the Toeplitz operator on  $H^{p} ,$  defined by

     $$T_{f} h=\int _{{\rm {\mathbb T}}} \frac{\overline{f}(\zeta )h(\zeta )}{1-\overline{\zeta }z} dm(\zeta ),{\kern 1pt} {\kern 1pt} {\kern 1pt} {\kern 1pt} {\kern 1pt} {\kern 1pt} {\kern 1pt} {\kern 1pt} {\kern 1pt} {\kern 1pt} {\kern 1pt} {\kern 1pt} {\kern 1pt} {\kern 1pt} {\kern 1pt} {\kern 1pt} {\kern 1pt} {\kern 1pt} {\kern 1pt} h\in H^{p} .$$

     In this paper we prove some sufficient conditions for the sequences of numbers   $\alpha =\left(\alpha _{n} \right)_{n\ge 0} $  in which the functions

  $$f*\alpha {\kern 1pt} {\kern 1pt} {\kern 1pt} {\kern 1pt} {\kern 1pt} {\kern 1pt} {\kern 1pt} {\kern 1pt} {\kern 1pt} {\kern 1pt} {\kern 1pt} {\kern 1pt} {\kern 1pt} {\kern 1pt} \mathop{=}\limits^{def} {\kern 1pt} {\kern 1pt} {\kern 1pt} {\kern 1pt} {\kern 1pt} {\kern 1pt} {\kern 1pt} {\kern 1pt} \sum \limits _{n\ge 1}\hat{f}(n){\kern 1pt} {\kern 1pt} \alpha _{n}  {\kern 1pt} {\kern 1pt} z^{n} {\kern 1pt} $$
determine bounded Toeplitz operators  $T_{f*\alpha } $ {\it  }on  $H^{1} $  and  $H^{\infty } $  for all  $f\in H^{\infty } .$
\end{abstract}
\section{Introduction}

     Let  $A$  be the class of all functions analytic in the unit disk   ${\rm {\mathbb D}}=\left\{\zeta :{\kern 1pt} {\kern 1pt} {\kern 1pt} {\kern 1pt} {\kern 1pt} {\kern 1pt} \left|\zeta \right|{\kern 1pt} <1\right\},$   $m(\zeta )$ - normalized Lebesgue measure on the circle  ${\rm {\mathbb T}}=\left\{\zeta :{\kern 1pt} {\kern 1pt} {\kern 1pt} {\kern 1pt} {\kern 1pt} {\kern 1pt} \left|\zeta \right|{\kern 1pt} =1\right\}.$ Let  $H^{p{\kern 1pt} {\kern 1pt} {\kern 1pt} {\kern 1pt} } (0<p\le \infty )$  is the  space of all functions analytic in  ${\rm {\mathbb D}}$  and satisfying
     $$\left\| f\right\| _{H^{p} }^{p} =\mathop{\sup }\limits_{0<{\kern 1pt} {\kern 1pt} r<{\kern 1pt} {\kern 1pt} {\kern 1pt} 1} \int _{{\rm {\mathbb T}}}\left|f(r\zeta )\right| ^{p} dm(\zeta )<\infty ,{\kern 1pt} {\kern 1pt} {\kern 1pt} {\kern 1pt} {\kern 1pt} {\kern 1pt} {\kern 1pt} {\kern 1pt} {\kern 1pt} {\kern 1pt} {\kern 1pt} {\kern 1pt} {\kern 1pt} {\kern 1pt} {\kern 1pt} {\kern 1pt} {\kern 1pt} {\kern 1pt} {\kern 1pt} {\kern 1pt} {\kern 1pt} {\kern 1pt} {\kern 1pt} {\kern 1pt} 0<p<\infty ,$$
     $$\left\| f\right\| _{H^{\infty } } =\mathop{\sup }\limits_{z\in {\rm {\mathbb D}}} {\kern 1pt} {\kern 1pt} {\kern 1pt} {\kern 1pt} \left|f(z)\right|<\infty ,{\kern 1pt} {\kern 1pt} {\kern 1pt} {\kern 1pt} {\kern 1pt} {\kern 1pt} {\kern 1pt} {\kern 1pt} {\kern 1pt} {\kern 1pt} {\kern 1pt} {\kern 1pt} {\kern 1pt} {\kern 1pt} {\kern 1pt} {\kern 1pt} {\kern 1pt} {\kern 1pt} {\kern 1pt} {\kern 1pt} {\kern 1pt} {\kern 1pt} {\kern 1pt} {\kern 1pt} p=\infty .$$

     Let  $M$  is the space of all finite, complex Borel measures on  ${\rm {\mathbb T}}$  with the usual variation norm.

     For $\mu \in M$ , the analytic function on  ${\rm {\mathbb D}}$

     $$K_{\mu } (z)=\int _{{\rm {\mathbb T}}}\frac{1}{1-\overline{\zeta }z}  {\kern 1pt} {\kern 1pt} {\kern 1pt} {\kern 1pt} {\kern 1pt} d\mu (\zeta )$$
is called the Cauchy transforms of  $\mu $  and the set of functions
$$K=\left\{f\in A:{\kern 1pt} {\kern 1pt} {\kern 1pt} {\kern 1pt} {\kern 1pt} {\kern 1pt} {\kern 1pt} {\kern 1pt} {\kern 1pt} f=K_{\mu } ,{\kern 1pt} {\kern 1pt} {\kern 1pt} {\kern 1pt} {\kern 1pt} {\kern 1pt} {\kern 1pt} {\kern 1pt} {\kern 1pt} {\kern 1pt} {\kern 1pt} \mu \in M\right\}$$
is called the space of Cauchy transforms.

     For  $d\mu (\zeta )=\varphi (\zeta )dm(\zeta ),{\kern 1pt} {\kern 1pt} {\kern 1pt} {\kern 1pt} {\kern 1pt} {\kern 1pt} {\kern 1pt} {\kern 1pt} {\kern 1pt} {\kern 1pt} {\kern 1pt} {\kern 1pt} \varphi \in L^{p} {\kern 1pt} ,{\kern 1pt} {\kern 1pt} {\kern 1pt} {\kern 1pt} {\kern 1pt} {\kern 1pt} {\kern 1pt} {\kern 1pt} {\kern 1pt} {\kern 1pt} 1\le p\le \infty $ , we denote  $K_{\mu } (z)=K_{\varphi } (z)$  and

     $$K^{p} =\left\{f\in A:{\kern 1pt} {\kern 1pt} {\kern 1pt} {\kern 1pt} {\kern 1pt} {\kern 1pt} {\kern 1pt} {\kern 1pt} {\kern 1pt} f=K_{\varphi } ,{\kern 1pt} {\kern 1pt} {\kern 1pt} {\kern 1pt} {\kern 1pt} {\kern 1pt} {\kern 1pt} {\kern 1pt} {\kern 1pt} {\kern 1pt} {\kern 1pt} \varphi \in L^{p} {\kern 1pt} \right\},{\kern 1pt} {\kern 1pt} {\kern 1pt} {\kern 1pt} {\kern 1pt} {\kern 1pt} {\kern 1pt} {\kern 1pt} 1\le p\le \infty .$$

     By the theorem of M. Riez  $K^{p} =H^{p} $  for ${\kern 1pt} {\kern 1pt} {\kern 1pt} {\kern 1pt} {\kern 1pt} 1<p<\infty $ , however  $H^{1} {\rm \varsubsetneq }K^{1} $ ,   $H^{\infty } {\rm \varsubsetneq }K^{\infty } $ .

     We note that  $K^{\infty } =BMOA$  (the space of analytic functions of bounded mean  oscillation )[1].

     For  $f\in H^{\infty } $  we denote by  $T_{f} $  the Toeplitz operator on  $H^{p} $ , defined by

     $$T_{f} h=K_{\overline{f}h} (z)=\int _{{\rm {\mathbb T}}} \frac{\overline{f}(\zeta )h(\zeta )}{1-\overline{\zeta }z} dm(\zeta ),{\kern 1pt} {\kern 1pt} {\kern 1pt} {\kern 1pt} {\kern 1pt} {\kern 1pt} {\kern 1pt} {\kern 1pt} {\kern 1pt} {\kern 1pt} {\kern 1pt} {\kern 1pt} {\kern 1pt} {\kern 1pt} {\kern 1pt} {\kern 1pt} {\kern 1pt} {\kern 1pt} {\kern 1pt} h\in H^{p} .$$

     By the theorem of M. Riez for ${\kern 1pt} {\kern 1pt} {\kern 1pt} {\kern 1pt} {\kern 1pt} 1<p<\infty $  the operator{\it  } $T_{f} $  is bounded on  $H^{p} $  for all  $f\in H^{\infty } .$  But if  ${\kern 1pt} {\kern 1pt} {\kern 1pt} {\kern 1pt} {\kern 1pt} p=1$  and  ${\kern 1pt} {\kern 1pt} {\kern 1pt} {\kern 1pt} {\kern 1pt} p=\infty $  not every function  $f\in H^{\infty } $  gives rise to bounded Toeplitz operator  $T_{f} $  on  $H^{1} $  and $H^{\infty } $ .

     There is also an interesting connection between multipliers of the spaces  $K$  and  $K^{p} $ , ${\kern 1pt} {\kern 1pt} {\kern 1pt} {\kern 1pt} {\kern 1pt} p=1,{\kern 1pt} {\kern 1pt} {\kern 1pt} {\kern 1pt} \infty $   and the Toeplitz operators.

     Let  ${\rm {\mathfrak M}}$  and  ${\rm {\mathfrak M}}^{p} $  be the class to all multipliers of the spaces  $K$  and  $K^{p} $ :
\
\
     $${\rm {\mathfrak M}}=\left\{f\in A:{\kern 1pt} {\kern 1pt} {\kern 1pt} {\kern 1pt} {\kern 1pt} {\kern 1pt} {\kern 1pt} {\kern 1pt} {\kern 1pt} f{\kern 1pt} g\in K,{\kern 1pt} {\kern 1pt} {\kern 1pt} {\kern 1pt} {\kern 1pt} {\kern 1pt} {\kern 1pt} {\kern 1pt} {\kern 1pt} {\kern 1pt} \forall {\kern 1pt} g\in K\right\},$$

     $${\rm {\mathfrak M}}^{p} =\left\{f\in A:{\kern 1pt} {\kern 1pt} {\kern 1pt} {\kern 1pt} {\kern 1pt} {\kern 1pt} {\kern 1pt} {\kern 1pt} {\kern 1pt} f{\kern 1pt} g\in K^{p} ,{\kern 1pt} {\kern 1pt} {\kern 1pt} {\kern 1pt} {\kern 1pt} {\kern 1pt} {\kern 1pt} {\kern 1pt} {\kern 1pt} {\kern 1pt} \forall {\kern 1pt} g\in K^{p} \right\}.$$

     Since  $K^{p} =H^{p} $  for ${\kern 1pt} {\kern 1pt} {\kern 1pt} {\kern 1pt} {\kern 1pt} 1<p<\infty $ , then  ${\rm {\mathfrak M}}^{p} =H^{\infty } $  for ${\kern 1pt} {\kern 1pt} {\kern 1pt} {\kern 1pt} {\kern 1pt} 1<p<\infty .$

     However  $${\rm {\mathfrak M}}={\rm {\mathfrak M}}^{1} {\rm \varsubsetneq }H^{\infty }  ,{\kern 1pt} {\kern 1pt} {\kern 1pt}{\kern 1pt} {\kern 1pt}{\kern 1pt} {\kern 1pt}{\kern 1pt} {\rm {\mathfrak M}}^{\infty } {\rm \subsetneqq }H^{\infty } $$  and

     $$\displaystyle {\rm {\mathfrak M}}={\rm {\mathfrak M}}^{1} =\left\{f\in H^{\infty } :{\kern 1pt} {\kern 1pt} {\kern 1pt}  {\kern 1pt} {\kern 1pt} {\kern 1pt} \left\| T_{f} \right\| _{H^{\infty } } <\infty \right\}[3],$$
         $$\displaystyle {\rm {\mathfrak M}}^{\infty } =\left\{f\in H^{\infty } :{\kern 1pt} {\kern 1pt} {\kern 1pt} {\kern 1pt} {\kern 1pt} {\kern 1pt} {\kern 1pt} {\kern 1pt} \left\| T_{f} \right\| _{H^{1} } <\infty \right\}[2].$$

     Let's note, that more information, bibliography and  review of results for  the spaces  $K$ and  ${\rm {\mathfrak M}}$  contains the  new monograph [5] .

     Since   ${\rm {\mathfrak M}}={\rm {\mathfrak M}}^{1} {\rm \varsubsetneq }H^{\infty } $ ,  ${\rm {\mathfrak M}}^{\infty } {\rm \varsubsetneq } H^{\infty } $  i.e. not all function   $f\in H^{\infty } $  give rise to bounded Toeplitz operators on  $H^{1} $  and  $H^{\infty } ,$  then naturally arises the following task:

     {\it To describe the}{\it  sequences}{\it  }{\it of }{\it  }{\it numbers}{\it  } $\alpha =\left(\alpha _{n} \right)_{n\ge 0} $ {\it , for which}{\it  the function}{\it s}

 $$f*\alpha {\kern 1pt} {\kern 1pt} {\kern 1pt} {\kern 1pt} {\kern 1pt} {\kern 1pt} {\kern 1pt} {\kern 1pt} {\kern 1pt} {\kern 1pt} {\kern 1pt} {\kern 1pt} {\kern 1pt} {\kern 1pt} \mathop{=}\limits^{def} {\kern 1pt} {\kern 1pt} {\kern 1pt} {\kern 1pt} {\kern 1pt} {\kern 1pt} {\kern 1pt} {\kern 1pt} \sum \limits _{n\ge 1}\hat{f}(n){\kern 1pt} {\kern 1pt} \alpha _{n}  {\kern 1pt} {\kern 1pt} z^{n} {\kern 1pt},{\kern 1pt}{\kern 1pt} {\kern 1pt} {\kern 1pt} {\kern 1pt} {\kern 1pt}  z\in {\rm {\mathbb D}}$$
\

\
{\it give }{\it  ri}{\it se to bounded Toeplitz operators}{\it  } $T_{f*\alpha } $ {\it  }{\it on } $H^{1} $ {\it  and } $H^{\infty } $ {\it  for all } $f\in H^{\infty } .$
\

\

     In this paper we prove some sufficient conditions for the sequences  $\alpha =\left(\alpha _{n} \right)_{n\ge 0} $  in which Toeplitz operator{\it  } $T_{f*\alpha } $  is bounded on  $H^{1} $  and  $H^{\infty } $  for all  $f\in H^{\infty } $ .

     Further we will use the following important theorem:
\

\

     {\bf Theorem of Smirnov.}{\bf  }

     {\it Let } ${\kern 1pt} {\kern 1pt} {\kern 1pt} {\kern 1pt} {\kern 1pt} 0<p<q$ {\it , } $f\in H^{p} $ {\it  and has } $L^{q} $ {\it  boundary values }{\it (} $f\in L^{q} ({\rm {\mathbb T}})$ {\it ). Then } $f\in H^{q} $ {\it .}
     \

     \

     We include also its proof for convenience of the reader.

     {\it Proof.} Since  $f\in H^{p} ,$ {\it  }then  $f=Bg$ {\it , }where  $B$ is a Blaschke product,{\it  } $g\in H^{p} $ and  $g\ne 0$  in  ${\rm {\mathbb D}}$ .

     The function  $g^{p} \in H^{1} $  and applying the formula of Poisson to the function  $g^{p} $ {\it  }we have

      $$g^{p} (z)=\int _{{\rm {\mathbb T}}}g^{p} (\zeta )P_{z} (\zeta ) {\kern 1pt} {\kern 1pt} {\kern 1pt} {\kern 1pt} {\kern 1pt} dm(\zeta ),{\kern 1pt} {\kern 1pt} {\kern 1pt} {\kern 1pt} {\kern 1pt}{\kern 1pt} {\kern 1pt} {\kern 1pt}{\kern 1pt} P_{z} (\zeta )=\frac{1-\left|z\right|^{2} }{\left|\zeta -z\right|^{2} } ,{\kern 1pt} {\kern 1pt} {\kern 1pt} {\kern 1pt} {\kern 1pt} {\kern 1pt} {\kern 1pt} {\kern 1pt} {\kern 1pt} {\kern 1pt} \zeta \in {\rm {\mathbb T}},{\kern 1pt} {\kern 1pt} {\kern 1pt} {\kern 1pt} {\kern 1pt} {\kern 1pt} {\kern 1pt} {\kern 1pt} {\kern 1pt} z\in D.$$

     From this formula, taking into account that

      $\left|f(z)\right|\le \left|g(z)\right|$ {\it  }{\it  }in  ${\rm {\mathbb D}}$ ,       $\left|f(\zeta )\right|=\left|g(\zeta )\right|$ {\it    }{\it  }for almost{\it  }every{\it  }{\it   } $\zeta \in {\rm {\mathbb T}}$ ,

follows

     $$\left|f(z)\right|^{p} \le \int _{{\rm {\mathbb T}}}\left|f(\zeta )\right|^{p} P_{z} (\zeta ) {\kern 1pt} {\kern 1pt} {\kern 1pt} {\kern 1pt} {\kern 1pt} dm(\zeta ).$$

     If ${\kern 1pt} {\kern 1pt} {\kern 1pt} {\kern 1pt} {\kern 1pt} q=\infty $ {\it ,} then    $f\in L^{\infty } ({\rm {\mathbb T}})$    and    $\left\| f\right\| _{H^{\infty } } \le \left\| f\right\| _{L^{\infty } ({\rm {\mathbb T}})} <\infty .$

     If ${\kern 1pt} {\kern 1pt} {\kern 1pt} {\kern 1pt} {\kern 1pt} q<\infty $ {\it , }then applying the Holder's inequality{\it  }we have

     $$\left|f(z)\right|^{p} \le \int _{{\rm {\mathbb T}}}\left|f(\zeta )\right|^{p} \left(P_{z} (\zeta )\right) {\kern 1pt} ^{p/q} {\kern 1pt} {\kern 1pt} {\kern 1pt} \left(P_{z} (\zeta )\right)^{1-p/q} {\kern 1pt} dm(\zeta )\le $$

     $$\le \left(\int _{{\rm {\mathbb T}}}{\kern 1pt} {\kern 1pt} \left|f(\zeta )\right|^{q}  {\kern 1pt} {\kern 1pt} {\kern 1pt} {\kern 1pt} P_{z} (\zeta )dm(\zeta )\right)^{p/q} \left(\int _{{\rm {\mathbb T}}}{\kern 1pt} {\kern 1pt}  {\kern 1pt} {\kern 1pt} P_{z} (\zeta ){\kern 1pt} dm(\zeta )\right)^{1-p/q} =$$

      $$=\left(\int _{{\rm {\mathbb T}}}{\kern 1pt} {\kern 1pt} \left|f(\zeta )\right|^{q}  {\kern 1pt} {\kern 1pt} {\kern 1pt} {\kern 1pt} P_{z} (\zeta )dm(\zeta )\right)^{p/q} {\kern 14pt} {\kern 4pt}\Rightarrow$$

    $$\left|f(z)\right|^{q} \le \int _{{\rm {\mathbb T}}}{\kern 1pt} {\kern 1pt} \left|f(\zeta )\right|^{q}  {\kern 1pt} {\kern 1pt} {\kern 1pt} {\kern 1pt} P_{z} (\zeta )dm(\zeta ).$$

     Integrating on the circle   $\left|z\right|{\kern 1pt} =r,{\kern 1pt} {\kern 1pt} {\kern 1pt} {\kern 1pt} {\kern 1pt} {\kern 1pt} {\kern 1pt} {\kern 1pt} 0<r<1$   we obtain

     $$\int _{{\rm {\mathbb T}}}\left|f(r\eta )\right|^{q}  dm(\eta )\le \int _{{\rm {\mathbb T}}}\int _{{\rm {\mathbb T}}}{\kern 1pt} {\kern 1pt} \left|f(\zeta )\right|^{q}  {\kern 1pt} {\kern 1pt} {\kern 1pt} {\kern 1pt} \frac{1-r^{2} }{\left|\zeta -r\eta \right|^{2} } dm(\zeta ) dm(\eta )\le \left\| f\right\| _{L^{q} ({\rm {\mathbb T}})} <\infty .$$

     Consequently  $f\in H^{q} .\Box$
\section{Main results}

     Let  ${\rm {\mathfrak N}}$  is the class of all functions  $f\in H^{\infty } $  for which

     $$\Lambda (f){\kern 1pt} {\kern 1pt} {\kern 1pt} {\kern 1pt} {\kern 1pt} \mathop{=}\limits^{def} {\kern 1pt} {\kern 1pt} {\kern 1pt} {\kern 1pt} {\kern 1pt} {\kern 1pt} \mathop{ess\sup }\limits_{\eta \in {\rm {\mathbb T}}} \int _{{\rm {\mathbb T}}} \frac{\left|f(\zeta )-f(\eta )\right|}{\left|\zeta -\eta \right|} dm(\zeta ){\kern 1pt} {\kern 1pt} {\kern 1pt} {\kern 1pt} {\kern 1pt} {\kern 1pt} {\kern 1pt} {\kern 1pt} <\infty {\kern 1pt} {\kern 1pt} {\kern 1pt} {\kern 1pt} .$$

For $f\in {\rm {\mathfrak N}}$  we denote   $\left\| f\right\| _{{\rm {\mathfrak N}}} {\kern 1pt} {\kern 1pt} {\kern 1pt} {\kern 1pt} \mathop{=}\limits^{def} {\kern 1pt} {\kern 1pt} {\kern 1pt} {\kern 1pt} {\kern 1pt} \left\| f\right\| _{H^{\infty } } +\Lambda (f).$
\

\

     {\bf Theorem 1.}  {\it If } $f\in {\rm {\mathfrak N}}{\kern 1pt} $ {\it , }{\it then Toeplitz operator } $T_{f} $ {\it  is bounded on } $H^{p} $  ${\kern 1pt} {\kern 1pt} {\kern 1pt} {\kern 1pt} ({\kern 1pt} p=1,{\kern 1pt} {\kern 1pt} {\kern 1pt} {\kern 1pt} {\kern 1pt} \infty )$ {\it  and}
     $${\kern 1pt} {\kern 1pt} \left\| T_{f} \right\| _{H^{p} } \le \left\| f\right\| _{{\rm {\mathfrak N}}} {\kern 1pt} {\kern 1pt} .$$

     {\it Proof.} The case  ${\kern 1pt} {\kern 1pt} {\kern 1pt} {\kern 1pt} {\kern 1pt} p=\infty $  is proved in [3,4] and is generalized in [6] for the multipliers of the integrals of Cauchy-Stieltjes type in domains with closed Jordan curve.

     We shall prove the case ${\kern 1pt} {\kern 1pt} {\kern 1pt} {\kern 1pt} {\kern 1pt} p=1$ .

     Let $f\in {\rm {\mathfrak N}}{\kern 1pt} $ {\it ,} ${\kern 1pt} {\kern 1pt} {\kern 1pt} {\kern 1pt} {\kern 1pt} h\in H^{1} $ . Let  ${\rm {\mathbb E}}$  be a subset with total measure  $\left(m\left({\rm {\mathbb E}}\right)=1\right)$  lying on  ${\rm {\mathbb T}}$  so that

$$\left\| f\right\| _{H^{\infty } } =\mathop{\sup }\limits_{\eta \in {\rm {\mathbb E}}{\kern 1pt} } {\kern 1pt} {\kern 1pt} {\kern 1pt} \left|f(\eta )\right|.$$

     Then

     $$\left\| T_{f} h\right\| _{H^{1} } =\mathop{\sup }\limits_{0<r<1} \int _{{\rm {\mathbb T}}}\left|\int _{{\rm {\mathbb T}}} \frac{\overline{f}(\zeta )h(\zeta )}{\zeta -r\eta } \zeta dm(\zeta )\right| dm(\eta )=$$

     $$=\mathop{\sup }\limits_{0<r<1} \int _{{\rm {\mathbb T}}}\left|\int _{{\rm {\mathbb T}}} \frac{\overline{f}(\zeta )-\overline{f}(r\eta )}{\zeta -r\eta } h(\zeta )\zeta dm(\zeta )+\overline{f}(r\eta )\int _{{\rm {\mathbb T}}} \frac{1}{\zeta -r\eta } h(\zeta )\zeta dm(\zeta )\right| dm(\eta )\le $$

     $$\le \mathop{\sup }\limits_{0<r<1} \left\{\int _{{\rm {\mathbb T}}}\int _{{\rm {\mathbb T}}} \left|\frac{\overline{f}(\zeta )-\overline{f}(r\eta )}{\zeta -r\eta } \right|\left|h(\zeta )\right|dm(\zeta ){\kern 1pt} {\kern 1pt} dm(\eta )+\int _{{\rm {\mathbb T}}}\left|\overline{f}(r\eta )h(r\eta )\right|  {\kern 1pt} dm(\eta )\right\}\le $$

     $$\le \mathop{\sup }\limits_{\begin{array}{l} {0<r<1} \\ {\zeta \in {\rm {\mathbb E}}} \end{array}} \left(\int _{{\rm {\mathbb T}}}{\kern 1pt} {\kern 1pt} \left|\frac{f(\zeta )-f(r\eta )}{\zeta -r\eta } \right|dm(\eta )+\left\| f\right\| _{H^{\infty } }  \right)\left\| h\right\| _{H^{1} } .$$

     We denote for  $\zeta \in {\rm {\mathbb E}}$

       $$F_{\zeta } (z)=\frac{f(\zeta )-f(z)}{\zeta -z} ,{\kern 1pt} {\kern 1pt} {\kern 1pt} {\kern 1pt} {\kern 1pt} {\kern 1pt} {\kern 1pt} {\kern 1pt} {\kern 1pt} {\kern 1pt} z\in {\rm {\mathbb D}}.$$

     Then

     $$\left\| T_{f} h\right\| _{H^{1} } \le \mathop{\sup }\limits_{\zeta \in {\rm {\mathbb E}}} \left(\left\| F_{\zeta } \right\| _{H^{1} } +\left\| f\right\| _{H^{\infty } } \right)\left\| h\right\| _{H^{1} } .$$

     To end the proof is necessary to show

     $$f\in {\rm {\mathfrak N}}{\kern 1pt} {\kern 1pt} {\kern 1pt} {\kern 1pt} {\kern 1pt} {\kern 1pt} {\kern 1pt} {\kern 1pt} \Rightarrow {\kern 1pt} {\kern 1pt} {\kern 1pt} {\kern 1pt} {\kern 1pt} \mathop{\sup }\limits_{\zeta \in {\rm {\mathbb E}}} \left\| F_{\zeta } \right\| _{H^{1} } <\infty .{\kern 1pt} $$

     Since

      $$\frac{1}{\zeta -z} \in H^{p} {\kern 1pt} {\kern 1pt} {\kern 1pt} {\kern 1pt} {\kern 1pt} (0<p<1)$$  and $f\in H^{\infty } $ ,  then   $F_{\zeta } (z)\in H^{p} {\kern 1pt} {\kern 1pt} {\kern 1pt} {\kern 1pt} {\kern 1pt} (0<p<1)$ .

     Furthermore

     $$f\in {\rm {\mathfrak N}}{\kern 1pt} {\kern 1pt} {\kern 1pt} {\kern 1pt} {\kern 1pt} {\kern 1pt} {\kern 1pt} {\kern 1pt} \Rightarrow {\kern 1pt} {\kern 1pt} {\kern 1pt} {\kern 1pt} {\kern 1pt} \mathop{\sup }\limits_{\zeta \in {\rm {\mathbb E}}} \left\| F_{\zeta } \right\| _{L^{1} ({\rm {\mathbb T}})} \le \Lambda (f)<\infty {\kern 1pt} {\kern 1pt} {\kern 1pt} {\kern 1pt} $$

and according to the Theorem of Smirnov
\

\

 $F_{\zeta } (z)\in H^{1} {\kern 1pt} {\kern 1pt} $ ,    $\left\| F_{\zeta } \right\| _{H^{1} } =\left\| F_{\zeta } \right\| _{L^{1} ({\rm {\mathbb T}})} \le \Lambda (f)<\infty $ .
\

\

     Consequently
\

\

     $$\left\| T_{f} \right\| _{H^{1} } \le \mathop{\sup }\limits_{\zeta \in {\rm {\mathbb E}}} \left(\left\| F_{\zeta } \right\| _{H^{1} } +\left\| f\right\| _{H^{\infty } } \right)\le \Lambda (f)+\left\| f\right\| _{H^{\infty } } =\left\| f\right\| _{{\rm {\mathfrak N}}} <\infty .\Box$$

     {\bf Remark. }{\it We note that from  the Theorem of Stegenga [2]  characterizing a class  of bounded Toeplitz operators on  $H^{1} $  does not follow Theorem 1 for ${\kern 1pt} {\kern 1pt} {\kern 1pt} {\kern 1pt} {\kern 1pt} p=1.$
\

\

     {\bf Lemma 1.}{\bf  [}{\bf 3}{\bf ] }{\it If } $p_{n} $ {\it  is a polynomial of degree} $n$ {\it , then}
$$\left\| p_{n} \right\| _{{\rm {\mathfrak N}}} {\kern 1pt} \le 3{\kern 1pt} {\kern 1pt} {\kern 1pt} {\kern 1pt} \left\| p_{n} \right\| _{H^{\infty } } \log (n+2).$$

     {\bf Definition.} {\it A sequence } $\alpha =\left(\alpha _{n} \right)_{n\ge 0} $ {\it  }{\it of positive numbers is called concave if}
      $$\alpha _{n+2} -\alpha _{n+1} \ge \alpha _{n+1} -\alpha _{n}  {\it    } \Leftrightarrow  {\it    } \alpha _{n} -2\alpha _{n+1} +\alpha _{n+2} \ge 0.$$

     {\bf Theorem 2}{\bf .} {\it Let } $\alpha =\left(\alpha _{n} \right)_{n\ge 0} $ {\it  be a monotone decreasing, concave sequence of positive numbers and  }
     $$\left\| \alpha {\kern 1pt} {\kern 1pt} \right\| {\kern 1pt} {\kern 1pt} {\kern 1pt} {\kern 1pt} {\kern 1pt} {\kern 1pt} {\kern 1pt} {\kern 1pt} {\kern 1pt} {\kern 1pt} {\kern 1pt} {\kern 1pt} \mathop{=}\limits^{def} {\kern 1pt} {\kern 1pt} {\kern 1pt} {\kern 1pt} {\kern 1pt} {\kern 1pt} {\kern 1pt} {\kern 1pt} \sum \limits _{n\ge 0}\frac{\alpha _{n} }{n+1} <\infty . {\kern 1pt} $$
     {\it Then} $f*\alpha {\kern 1pt} {\kern 1pt} {\kern 1pt} \in {\rm {\mathfrak N}}$ {\it , Toeplitz operator}{\it  } $T_{f*\alpha } $ {\it  is bounded on } $H^{1} $ {\it  and } $H^{\infty } $ {\it  for all } $f\in H^{\infty } $ {\it and}
$$\left\| T_{f*\alpha } \right\| _{H^{p} } \le \left\| f*\alpha \right\| _{{\rm {\mathfrak N}}} \le 12\left\| f\right\| _{H^{\infty } } \left\| \alpha {\kern 1pt} {\kern 1pt} \right\| {\kern 1pt} {\kern 1pt} {\kern 1pt} ,{\kern 1pt} {\kern 1pt} {\kern 1pt} {\kern 1pt} {\kern 1pt} {\kern 1pt} {\kern 1pt} {\kern 1pt} {\kern 1pt} p=1,\infty .$$}

     {\it Proof. } Using Abel's formula two times we obtain

     $${\kern 1pt} {\kern 1pt} {\kern 1pt} \sum \limits _{n\ge 0}\frac{\alpha _{n} }{n+1} =\sum \limits _{n\ge 0}(\alpha _{n} -\alpha _{n+1} )\sum \limits _{k=0}^{n}\frac{1}{k+1}    {\kern 1pt} {\kern 1pt} {\kern 1pt} {\kern 1pt} {\kern 1pt} {\kern 1pt} \ge $$

$${\kern 1pt} \ge \sum \limits _{n\ge 0}(\alpha _{n} -\alpha _{n+1} )\log (n+2) {\kern 1pt} {\kern 1pt} {\kern 1pt} {\kern 1pt} {\kern 1pt} {\kern 1pt} =$$

$$=\sum \limits _{n\ge 0}(\alpha _{n} -2\alpha _{n+1} +\alpha _{n+2} )\sum \limits _{k=0}^{n}\log (k+2)  .$$

     Since

$$\sum \limits _{k=0}^{n}\log (k+2) \ge \sum \limits _{k=\left[n/2\right]}^{n}\log (k+2) \ge (n/2+1)\log (\left[n/2\right]+2)\ge \frac{1}{4} (n+1)\log (n+2),$$

 then

     $${\kern 1pt} {\kern 1pt} {\kern 1pt} 4\sum \limits _{n\ge 0}\frac{\alpha _{n} }{n+1}  {\kern 1pt} {\kern 1pt} {\kern 1pt} {\kern 1pt} {\kern 1pt} \ge \sum \limits _{n\ge 0}(\alpha _{n} -2\alpha _{n+1} +\alpha _{n+2} )(n+1)\log (n+2). $$

     Further let  $f\in H^{\infty } $  and

     $$\displaystyle {\kern 1pt} {\kern 1pt} {\kern 1pt} S_{n} (f)=\sum \limits _{k=0}^{n}\hat{f}(k){\kern 1pt} {\kern 1pt} z^{k}  ;{\kern 1pt} {\kern 1pt} {\kern 1pt} {\kern 1pt} {\kern 1pt} {\kern 1pt} {\kern 1pt} {\kern 1pt} {\kern 1pt} {\kern 1pt} {\kern 1pt} {\kern 1pt} {\kern 1pt} {\kern 1pt} {\kern 1pt} {\kern 1pt} {\kern 1pt} {\kern 1pt} {\kern 1pt} {\kern 1pt} {\kern 1pt} {\kern 1pt} {\kern 1pt} {\kern 1pt} {\kern 1pt} {\kern 1pt} {\kern 1pt} {\kern 1pt} {\kern 1pt} {\kern 1pt} \sigma _{n} (f)=\frac{1}{n+1} \sum \limits _{k=0}^{n}S_{k} (f) .$$

     Applying the Abel's formula we obtain

     $$f*\alpha {\kern 1pt} {\kern 1pt} {\kern 1pt} {\kern 1pt} {\kern 1pt} {\kern 1pt} {\kern 1pt} {\kern 1pt} {\kern 1pt} {\kern 1pt} {\kern 1pt} ={\kern 1pt} {\kern 1pt} {\kern 1pt} {\kern 1pt} {\kern 1pt} {\kern 1pt} {\kern 1pt} {\kern 1pt} \sum \limits _{n\ge 0}\hat{f}(n){\kern 1pt} {\kern 1pt} \alpha _{n}  {\kern 1pt} {\kern 1pt} z^{n} {\kern 1pt} {\kern 1pt} ={\kern 1pt} {\kern 1pt} {\kern 1pt} {\kern 1pt} {\kern 1pt} {\kern 1pt} {\kern 1pt} {\kern 1pt} \sum \limits _{n\ge 0}({\kern 1pt} {\kern 1pt} \alpha _{n}  -{\kern 1pt} \alpha _{n+1} ){\kern 1pt} {\kern 1pt} S_{n} (f)=$$

           $${\kern 1pt} =\sum \limits _{n\ge 0}(\alpha _{n} -2\alpha _{n+1} +\alpha _{n+2} )(n+1){\kern 1pt} {\kern 1pt} {\kern 1pt} \sigma _{n} (f). $$

     Since by Lemma 1.

     $$\left\| \sigma _{n} (f)\right\| _{{\rm {\mathfrak N}}} \le 3\left\| \sigma _{n} (f)\right\| _{H^{\infty } } \log (n+2)\le 3\left\| f\right\| _{H^{\infty } } \log (n+2),$$

     then

$$\left\| f*\alpha {\kern 1pt} \right\| {\kern 1pt} _{{\rm {\mathfrak N}}} {\kern 1pt} {\kern 1pt} {\kern 1pt} {\kern 1pt} {\kern 1pt} \le \sum \limits _{n\ge 0}(\alpha _{n} -2\alpha _{n+1} +\alpha _{n+2} )(n+1){\kern 1pt} {\kern 1pt} {\kern 1pt} \left\| \sigma _{n} (f)\right\| _{{\rm {\mathfrak N}}}  \le $$

$${\kern 1pt} {\kern 1pt} {\kern 1pt} \le 3\left\| f\right\| _{H^{\infty } } \sum \limits _{n\ge 0}(\alpha _{n} -2\alpha _{n+1} +\alpha _{n+2} )(n+1)\log (n+2) \le $$

$${\kern 1pt} {\kern 1pt} {\kern 1pt} \le 12\left\| f\right\| _{H^{\infty } } \sum \limits _{n\ge 0}\frac{\alpha _{n} }{n+1}  =12\left\| f\right\| _{H^{\infty } } \left\| \alpha \right\| <\infty .\Box$$

     The following proposition follows at once from Theorem 2.
     \

     \

     {\bf Theorem 3.} {\it Let } $\alpha $ {\it  denote one of the sequences  } $(\varepsilon >0)$ {\it  :}

     $$\displaystyle \left(\frac{1}{(n+1)^{\varepsilon } } \right)_{n\ge 0} ;{\kern 45pt}$$
      \
      \

     $$\displaystyle \left(\frac{1}{\log ^{1+\varepsilon } (n+2)} \right)_{n\ge 0} ;{\kern 45pt}$$
     \
     \

     $$\displaystyle \left(\frac{1}{\log (n+2)\log ^{1+\varepsilon } \log (n+3)} \right)_{n\ge 0} ,   ..............................$$
\

\

{\it Then} $f*\alpha {\kern 1pt} {\kern 1pt} {\kern 1pt} \in {\rm {\mathfrak N}}$ {\it , Toeplitz operator}{\it  } $T_{f*\alpha } $ {\it  is bounded on } $H^{1} $ {\it  and } $H^{\infty } $ {\it  for all } $f\in H^{\infty } .$

\

\
     {\bf Remark. }{\bf  }{\it Theorem 3  }{\it  was  proved by another method in }[3] ( {\it Theorem 7. ) for the bounded Toeplitz operators}{\it  } $T_{f*\alpha } $ {\it  on } $H^{\infty } .$
\

\

     {\bf Theorem 4}{\bf .} {\it Let the sequence } $\alpha =\left(\alpha _{n} \right)_{n\ge 0} $ {\it  satisfy}{\it  {\kern 4pt}the conditions of Theorem }{\it 3. If the sequence} $a=\left(a_{n} \right)_{n\ge 0} \in \ell ^{2} $ {\it , then there exist}{\it s}{\it  {\kern 4pt}a function} $f\in {\rm {\mathfrak N}}$ {\it , satisfying}

 $$\left|\hat{f}(n)\right|\ge \alpha _{n} \left|a_{n} \right|  , {\kern 25pt}     \left\| f\right\| _{{\rm {\mathfrak N}}} {\kern 1pt} \le c_{0} {\kern 1pt} {\kern 1pt} \left\| \alpha \right\| {\kern 1pt} {\kern 1pt} {\kern 1pt} \left\| a\right\| _{\ell ^{2} } ,$$

{\it where } $c_{0} {\kern 1pt} {\kern 1pt} $ {\it  is an absolute constant.}
\

\

     {\it Proof.} By the Theorem of Kislyakov [7] if $a=\left(a_{n} \right)_{n\ge 0} \in \ell ^{2} $ {\it ,} then there exists a function  $f\in H^{\infty } ,$  satisfying

 $$\left|\hat{g}(n)\right|\ge \left|a_{n} \right|  ,{\kern 25pt}\left\| g\right\| _{H^{\infty } } {\kern 1pt} \le B{\kern 1pt} {\kern 1pt} {\kern 1pt} {\kern 1pt} {\kern 1pt} \left\| a\right\| _{\ell ^{2} } ,$$

where  $B$  is an absolute constant. By Theorem 2.3  $f=g*\alpha {\kern 1pt} {\kern 1pt} {\kern 1pt} \in {\rm {\mathfrak N}}$ {\it  }and

$$\left\| f\right\| _{{\rm {\mathfrak N}}} \le 12\left\| \alpha \right\| \left\| g\right\| _{H^{\infty } } {\kern 1pt} \le 12B{\kern 1pt} {\kern 1pt} \left\| \alpha \right\| {\kern 1pt} {\kern 1pt} {\kern 1pt} \left\| a\right\| _{\ell ^{2} } .\Box$$

\

\

\

\

\

\

\noindent
{\small Department of Mathematics\\
        Technical University at Plovdiv\\
        25, Tsanko Dyustabanov,\\
        Plovdiv, Bulgaria\\
        e-mail: peyyyo@mail.bg}

.\\[4pt]

\begin{thebibliography}{199}


\bibitem{1} J. B. Garnett. {\it Bounded analytic functions.} Academic Press, Inc., New York-London, 1981. MR0628971 (83g:30037)

\bibitem{2} D. A. Stegenga. {\it Bounded Toeplitz operators on}  $H^{1} $ {\it  and applications of the duality between } $H^{1} $ {\it and the functions of bounded mean oscillation.} Amer. J. Math., 98, 1976,no. 3, 573-589. MR0420326 (54 \#8340)

\bibitem{3} S. A. Vinogradov. {\it Properties of multipliers of Cauchy - Stieltjes integrals and some factorization problems for analytic functions.} Amer. Math. Sos. Transl. (2) vol. 115, 1980, 1-32. MR0586560 (58 \#28518)

\bibitem{4}S. V. Hruscev, S. A. Vinogradov. {\it Inner functions and multipliers of Cauchy type integrals}. Ark. mat, 19, 1981, 23-42. MR0625535 (83c:30027)

\bibitem{5} J. A. Cima, A. L. Matheson, T. W. Ross. {\it The Cauchy transform.} American Mathematical Society, Providence, RI,{\it  }2006.{\it  }MR2215991 (2006m: 30003)

\bibitem{6} P. Stoilov. {\it Multipliers of integrals of Cauchy - Stieltjes type.} Mathematics and mathematical education, Publ. House Bulgar. Acad. Sci., Sofia, 1986, 316 - 319. (Russian) MR0872936 (88e:30104)

\bibitem{7} S. V. Kislyakov.{\bf  }{\it Fourier coefficients of boundary values of functions that are analytic in the disc and bidisc.}{\it  }Trudy Mat. Inst. Steklov. vol. 155,{\bf  }1981, 77--94. (Russian) MR0615566 (83a:42005)

\end{thebibliography}
\end{document}